\documentclass[a4paper,11pt,twoside]{article}
\usepackage{pst-fill,pst-grad}
\usepackage{textcomp}
\usepackage[english]{babel}
\usepackage{graphicx}
\usepackage{amsmath}
\usepackage{float}
\usepackage[matrix,arrow,curve]{xy}
\usepackage{pstricks}
\usepackage{amsmath,amsfonts,verbatim,afterpage,theorem,euscript,mathrsfs,amssymb}
\usepackage{amsfonts}
\usepackage{amssymb}
\usepackage{array}
\usepackage{dsfont}
\usepackage[colorlinks,citecolor=blue,hypertexnames=false,urlcolor=blue, linkcolor=blue, breaklinks=true]{hyperref}
\newtheorem{Theoreme}{Theorem}

\newtheorem{Proposition}{Proposition}[section]
\newtheorem{Lemme}{Lemma}[section]
\newtheorem{Corollaire}{Corollary}[section]
\newtheorem{Remarque}{Remark}[section]
\newcommand{\mysection}{\setcounter{equation}{0} \section}

\def\vu{\vec{u}}
\def\vv{\vec{v}}

\def\vpsi{\vec{\psi}}
\def\vn{\vec{\nabla}}
\def\R{\mathbb{R}^3}
\def\vf{\vec{f}}

\newcommand{\dv}{\operatorname{div}}

\newcommand\blfootnote[1]{%
	\begingroup
	\renewcommand\thefootnote{}\footnote{#1}%
	\addtocounter{footnote}{-1}%
	\endgroup
}

\title{\bf On an almost sharp Liouville type theorem for fractional Navier-Stokes equations} 
\author{Diego Chamorro\footnote{Laboratoire de Math\'ematiques et Mod\'elisation d'Evry, Universit\'e Paris-Saclay, France. \emph{diego.chamorro@univ-evry.fr}}, Bruno Poggi\footnote{Departament de  Matem\`atiques, Universitat Aut\`onoma de Barcelona, Bellaterra, Catalonia. \emph{poggi@mat.uab.cat}}}

\begin{document}
\maketitle
\begin{scriptsize}
\abstract{We investigate existence, Liouville type theorems and regularity results for the 3D stationary and incompressible fractional Navier-Stokes equations: in this setting the usual Laplacian is replaced by its fractional power $(-\Delta)^{\frac{\alpha}{2}}$ with $0<\alpha<2$. By applying a fixed point argument, weak solutions can be obtained in the Sobolev space $\dot{H}^{\frac{\alpha}{2}}(\R)$ and if we add an extra integrability condition, stated in terms of Lebesgue spaces, then we can prove for some values of $\alpha$ that the zero function is the unique smooth solution. The additional integrability condition is almost sharp for $3/5<\alpha<5/3$. Moreover, in the case $1<\alpha<2$ a gain of regularity is established under some conditions, however the study of regularity in the regime $0<\alpha\leq 1$ seems for the moment to be an open problem.}\\[3mm]
\end{scriptsize}
\blfootnote{The second author was supported by the European Research Council (ERC) under the European Union's Horizon 2020 research and innovation programme (grant agreement 101018680).}
\mysection{Introduction and presentation of the results}
In this article we   study existence, regularity and uniqueness properties of the 3D fractional Navier-Stokes equations which are given by the following system:
\begin{equation}\label{NSStationnaire1}
\begin{cases}
(-\Delta)^{\frac \alpha 2} \vu(x) +(\vu\cdot \vn)\vu(x)+\vn p(x)-\vf(x)=0,\quad\mbox{with } 0<\alpha<2,\\[3mm]
\dv(\vu)(x)=0, \quad x\in \R.
\end{cases}
\end{equation}
Here, the fractional operator $(-\Delta)^{\frac \alpha 2}$ is defined at the Fourier level by the symbol $|\xi|^\alpha$. Using the traditional notation, the vector field $\vu:\R\longrightarrow \R$ represents the velocity of the fluid, $p:\R\longrightarrow\mathbb{R}$ is the internal pressure of the fluid and $\vf:\R\longrightarrow \R$ is a given external force.\\

Before presenting our results related to the system (\ref{NSStationnaire1}), it is convenient to recall some facts about the usual stationary Navier-Stokes equations. Indeed, note that when $\alpha=2$, (\ref{NSStationnaire1}) is exactly the problem given by the classical incompressible Navier-Stokes equations
\begin{equation}\label{eq.nsc}
\begin{cases}
-\Delta\vu(x) +(\vu\cdot \vn)\vu(x)+\vn p(x)-\vf(x)=0,\\[3mm]
\dv(\vu)(x)=0, \quad x\in \R.
\end{cases}
\end{equation}
The problem (\ref{eq.nsc}) can be studied from different points of view; we first remark that the pressure $p$ can be easily deduced from the velocity field $\vu$ and the external force $\vf$ since, due to the divergence-free property of $\vu$, we have that
$$p=\frac{1}{(-\Delta)}\dv\big((\vu\cdot \vn)\vu-\vf\big),$$
and this fact allows us to focus our study on the velocity field $\vu$ (note that the same identity can be easily deduced from the system (\ref{NSStationnaire1}) since in both cases we have $\dv(\vu)=0$). Now, concerning existence problems for the usual Navier-Stokes equations (\ref{eq.nsc}), if we assume that $\vf\in \dot{H}^{-1}(\R)$ and that $\dv(\vf)=0$, then it is an easy exercise to construct solutions $\vu\in \dot{H}^1(\R)$ (see, for instance, \cite[Theorem 16.2]{PGLR2}) and moreover it is not hard to prove that these solutions are regular. However, a priori it is not known whether these solutions are unique, and an interesting open problem (initially mentioned in \cite{Galdi} and also stated in \cite{Ser2016}) is the following: show that any solution $\vu$ of the problem
\begin{equation}\label{NSStationnaire13}
-\Delta\vu+(\vu\cdot \vn)\vu+\vn p=0,
\end{equation}
which satisfies the conditions
\begin{equation}\label{Condition1}
\vu \in \dot{H}^1(\R) \qquad\mbox{and}\qquad \vu(x)\to 0 \mbox{ as } |x|\to +\infty,
\end{equation}
is identically equal to zero.\\ 

Note that, by the classical Sobolev embeddings, we have $\dot{H}^1(\R)\subset L^6(\R)$, but this seems not enough to conclude that a solution $\vu\in \dot{H}^1(\R)$ of the equation (\ref{NSStationnaire13}) is null.  Nevertheless, if we assume some additional hypotheses, for example $\vu\in E(\R)$ where $E$ is a nice functional space, then statements of the following form have been shown:
\begin{center}
if $\vu\in  \dot{H}^1(\R) \cap E(\R)$ is a solution of the equation (\ref{NSStationnaire13}) in $\mathbb R^3$, then we have $\vu\equiv 0$,
\end{center}
and this sort of result is known in the literature as a \emph{Liouville theorem} for the Navier-Stokes equations. In \cite{Galdi} the case $E=L^{\frac{9}{2}}(\R)$ was studied. The space $E=BMO^{-1}(\R)$ was considered in \cite{Koch} and other funcional spaces can also be taken into account, see for example the articles \cite{OJ1}, \cite{Kozono}, \cite{Ser2015}.\\

Remark that if we want to consider only one ``simple'' additional hypothesis, then a general Liouville-type theorem was proven in \cite{ChJL} with  $E=L^q(\R)$ for some 
\begin{equation}\label{ConditionLaplacienLiouville}
3\leq q\leq \frac{9}{2}.
\end{equation}
It is very interesting to note here that there is a \emph{gap} between this set of values and the integrability condition given in (\ref{Condition1}) -which is $\vu \in L^6(\R)$ due to the Sobolev embedding- as at present we do not know how to fill the distance between $\frac{9}{2}$ and 6. Thus the following problem: ``\emph{show that any solution $\vu$  of (\ref{NSStationnaire13}) with $\vu\in \dot{H}^1(\R)$ and $\vu\in L^q(\R)$ for some $\frac{9}{2}<q<6$, is identically equal to zero}'' remains, to the best of our knowledge, an open problem.\\

Let us come back now to the fractional Navier-Stokes equations (\ref{NSStationnaire1}). In particular, we are interested in understanding how the previous uniqueness results vary if we replace the Laplacian $\Delta$ by the operator $(-\Delta)^\frac{\alpha}{2}$ with $0<\alpha<2$. In \cite{wx18}, the authors use the Caffarelli-Silvestre extension \cite{cs07} to show that for $0<\alpha<2$, a  smooth weak solution $u\in\dot H^{\frac{\alpha}{2}}(\R)$ to (\ref{NSStationnaire1}) is trivial if $u\in L^{\frac92}(\mathbb R^3)$. On the other hand,  since $\dot H^{\frac{\alpha}{2}}(\mathbb R^3)$ embeds into $L^{\frac6{3-\alpha}}(\R)$, then for $\alpha<5/3$ it is reasonable to expect that the assumption $\vu\in L^{\frac92}(\mathbb R^3)$ may be  replaced by a more natural Lebesgue space whose exponent depends on the value of $\alpha$.\\


We begin by proving that, under some mild assumptions over the external force $\vf$, there exists at least one solution $\vu\in \dot{H}^{\frac{\alpha}{2}}(\R)$. Indeed we have:  
\begin{Theoreme}[Existence]\label{Theo_ExistenceStationnaires}
Fix $0< \alpha< 2$ and consider $\vf \in \dot{H}^{-1}(\R)\cap\dot{H}^{-\frac\alpha 2}(\R)$ an external force such that $\dv(\vf)=0$. There exists a divergence-free vector field $\vu\in  \dot{H}^{\frac \alpha 2}(\R)$ and a pressure $p\in \dot{H}^{\alpha-\frac32}(\R)$, such that $(\vu,p)$ is a solution of the stationary fractional Navier-Stokes equations (\ref{NSStationnaire1}).
\end{Theoreme}
Existence of certain weak solutions to the fractional Navier-Stokes equations (\ref{NSStationnaire1})  has already been studied in \cite{Tang} via the Caffarelli-Silvestre extension \cite{cs07}; our approach is quite different. We use the Schaefer fixed point theorem, which is a  useful tool when dealing with the existence of solutions for partial differential equations. In order to apply this general fixed point theorem, we will regularize the equation (\ref{NSStationnaire1}), and to recover the initial equation we will need to study a limit by considering  subsequences. This will give us a solution but we will lose uniqueness.\\

The study of the potential uniqueness of such solutions is in general a completely different open problem (besides the case $\alpha= 5/3$ which was studied in \cite{wx18}). However, if we add some extra conditions we can obtain interesting conclusions and in this sense we have our next result:
\begin{Theoreme}[Liouville type]\label{Theo_Liouville}
Consider the stationary fractional Navier-Stokes equations 
\begin{equation}\label{NSStationnaire2}
(-\Delta)^{\frac \alpha 2} \vu +(\vu\cdot \vn)\vu+\vn p=0,\qquad \dv(\vu)=0,\quad\qquad 0<\alpha<2.
\end{equation}
Assume that $\vu,p$ are smooth functions  that satisfy (\ref{NSStationnaire2}) and consider a positive parameter $0<\epsilon<2\alpha$.\\
\begin{itemize}
\item[1)] Let $\alpha=1$. If $\vu\in \dot{H}^{\frac1 2}(\R)\cap L^{\frac{6-\epsilon}{2}}(\R)$, then we have that $\vu=0$.\\[2mm]

\item[2)] Let $1< \alpha<2$ and fix the parameter $0<\epsilon<2\alpha$ such that
\begin{equation}\label{ConditionIntegrable}
1+\frac{\epsilon}{3}\leq \alpha\leq \frac{5}{3}+\frac{2}{9}\epsilon.
\end{equation}
If $\vu\in \dot{H}^{\frac\alpha 2}(\R)\cap L^{\frac{6-\epsilon}{3-\alpha}}(\R)$, then we have $\vu=0$.\\[2mm]
\item[3)] Let $\frac{3}{5}<\alpha<1$ and consider a parameter $0<\epsilon<2\alpha$ such that  
\begin{equation}\label{ConditionIntegrable1}
1-\frac{\epsilon}{3}\leq \alpha\leq \frac{5}{3}-\frac{2}{9}\epsilon.
\end{equation}
If $\vu\in \dot{H}^{\frac\alpha 2}(\R)\cap L^{\frac{6-\epsilon}{3-\alpha}}(\R)\cap L^{\frac{6+\epsilon}{3-\alpha}}(\R)$, then we have $\vu=0$.
\end{itemize}
\end{Theoreme}
Some remarks are in order. Indeed, we first note that the general condition $\vu\in \dot{H}^{\frac\alpha 2}(\R)$, stated in all the items above, is rather natural since from Theorem \ref{Theo_ExistenceStationnaires} we know how to construct solutions in this functional space. Second, by the classical Sobolev embeddings we have $\dot{H}^{\frac{\alpha}2}(\R)\subset  L^{\frac{6}{3-\alpha}}(\R)$ but we do not have $\vu \in L^{\frac{6-\epsilon}{3-\alpha}}(\R)$ nor $\vu \in L^{\frac{6+\varepsilon}{3-\alpha}}(\R)$ for $\epsilon>0$, and we can thus see that the conditions stated in the theorem are actual \emph{additional} hypotheses which help us to obtain this Liouville-type result. Next we note that if $\alpha\to 2$ then by the condition (\ref{ConditionIntegrable}) we have $\frac{3}{2}\leq \epsilon\leq 3$ and this lead us to the Lebesgue spaces $L^q(\R)$ with $3\leq q\leq \frac{9}{2}$, which is exactly the condition (\ref{ConditionLaplacienLiouville}) stated above and we recover the known results for the classical stationary Navier-Stokes equations regarding additional Lebesgue space hypotheses. We remark also that in the range $1\leq \alpha\leq \frac{5}{3}$ then, following the relationship (\ref{ConditionIntegrable}) (or (\ref{ConditionIntegrable1})), we can consider very small values for the parameter $\epsilon>0$ and thus the additional information $L^{\frac{6-\epsilon}{3-\alpha}}(\R)$ (or  $L^{\frac{6+\epsilon}{3-\alpha}}(\R)$) becomes closer and closer to the critical space $L^{\frac{6}{3-\alpha}}(\R)$: in the case of the stationary fractional Navier-Stokes equation we can almost fill the gap between the space $L^{\frac{6}{3-\alpha}}(\R)$ and the additional information required to deduce Liouville-type theorems. However, we can not simply take $\epsilon\to0$ as the information conveyed by the hypotheses (with $\epsilon>0$) is needed to obtain our results. Note finally that the lower limit $\frac{3}{5}$ stated in the third item is related to some technical issues. To finish, let us mention that we do not claim any optimality on the different relationships stated here.\\ 

To continue, we remark now that smoothness was taken for granted in the previous theorem, but this condition is redundant in some cases. Indeed, if we study the regularity of the solutions obtained in Theorem \ref{Theo_ExistenceStationnaires}, we have the following result.
\begin{Theoreme}[Regularity]\label{Theoreme_Regularite}
Consider the stationary fractional Navier-Stokes equations  (\ref{NSStationnaire2}).
\begin{itemize}
\item[1)] If $\frac{5}{3}<\alpha<2$, then the solutions  $\vu\in \dot{H}^{\frac{\alpha}{2}}(\R)$ obtained via Theorem \ref{Theo_ExistenceStationnaires} above are smooth. 

\item[2)] Let $1<\alpha\leq \frac{5}{3}$ and consider a solution $\vu\in \dot{H}^{\frac{\alpha}{2}}(\R)$. If we assume that $\vu\in L ^\infty(\R)$, then these solutions are smooth. 
\end{itemize}
\end{Theoreme}
We can see that the smoothness hypothesis in Theorem \ref{Theo_Liouville} is actually not necessary in the case $\frac{5}{3}<\alpha<2$. Nevertheless, if $1<\alpha\leq \frac{5}{3}$ the regularizing effect of the operator $(-\Delta)^{\frac{\alpha}{2}}$ seems to be too weak in order to obtain a gain of regularity and an additional hypothesis is thus warranted. For the sake of simplicity we assumed here a very strong condition, namely $\vu\in L ^\infty(\R)$, but we believe that other more general conditions can be considered. The study of the regularity in the case $0<\alpha\leq 1$ is considerably more difficult and technical to handle and, to the best of our knowledge, it constitutes an open problem that will not be treated here.\\

As a final remark, we point out that other functional spaces (such as Besov, Triebel-Lizorkin, Lorentz, Morrey spaces, etc.) can be used to develop all the previous theorems. However, the $L^2$-based Sobolev spaces are enough to highlight the behavior of the fractional Navier-Stokes equations considered here.\\ 

The plan of the article is the following: in Section \ref{Secc_Notation} we recall some notation and useful results. In Section \ref{Secc_ProofTheo1} we prove Theorem \ref{Theo_ExistenceStationnaires} and in Section \ref{Secc_ProofTheo2} we prove Theorem \ref{Theo_Liouville}. The last section is devoted to prove Theorem \ref{Theoreme_Regularite}.
 
\mysection{Preliminaries}\label{Secc_Notation}
For $1<p<+\infty$ and for $s>0$ we define the homogeneous Sobolev spaces $\dot{W}^{s,p}(\R)$ by the condition 
$$\|f\|_{\dot{W}^{s,p}}=\|(-\Delta)^{\frac{s}{2}}f\|_{L^p}<+\infty.$$
In the special case when $p=2$ we simply write $\dot{W}^{s,2}(\R)=\dot{H}^{s}(\R)$. The non-homogeneous Sobolev spaces $W^{s,p}(\R)$ are defined by the condition 
$$\|f\|_{W^{s,p}}=\|f\|_{L^p}+\|(-\Delta)^{\frac{s}{2}}f\|_{L^p}<+\infty,$$
from which we easily deduce the embedding $W^{s,p}(\R)\subset\dot{W}^{s,p}(\R)$. Note also that, if $s_1>s_0>0$ then we have the space inclusion $W^{s_1,p}(\R)\subset W^{s_0,p}(\R)$. As the Sobolev spaces will constitute our main framework, we recall in the following lemmas some classical and useful results. 
\begin{Lemme}[Sobolev embeddings]
\begin{itemize}
\item[]
\item[1)]For $0<s<\frac3p$ and $1<p,q<+\infty$, if we have the relationship $-\frac3q=s-\frac3p$, then   we have the classical Sobolev inequality
$$\|f\|_{L^q}\leq C\|f\|_{\dot{W}^{s,p}},\qquad\text{for each }f\in C_c^{\infty}(\mathbb R^n).$$
\item[2)]If $0<s_0<s_1$ and $1<p_0, p_1<+\infty$ are such that $s_0-\frac{3}{p_0}=s_1-\frac{3}{p_1}$, then we have the following Sobolev space inclusion:
$$\dot{W}^{s_1,p_1}(\R)\subset\dot{W}^{s_0,p_0}(\R).$$
\end{itemize}
\end{Lemme}
\begin{Lemme}[Rellich-Kondrachov]\label{Lemma_Rellich-Kondrachov}
Let $\Omega\subset \R$ be a bounded Lipschitz domain. If $0<s<\frac{3}{p}$, then for all $1\leq q<\frac{3p}{3-sp}$ we have the following compact inclusion
$$\dot{W}^{s,p}(\Omega)\subset\subset L^q(\Omega).$$
\end{Lemme}
A useful consequence of this lemma is that any uniformly bounded sequence in $\dot{W}^{s,p}(\Omega)$ has a subsequence that converges in $L^q(\Omega)$. 
\begin{Lemme}[Fractional Leibniz rule]\label{FracLeibniz}
\begin{itemize}
\item[]
\item[1)] Consider $f,g$ two smooth functions. Then we have the estimate
$$\|(-\Delta)^{\frac{s}{2}}(fg)\|_{L^p}\leq C\|(-\Delta)^{\frac{s}{2}}f\|_{L^{p_0}}\|g\|_{L^{p_1}}+C\|f\|_{L^{q_0}}\|(-\Delta)^{\frac{s}{2}}g\|_{L^{q_1}},$$
where $\frac1p=\frac{1}{p_0}+\frac{1}{p_1}=\frac{1}{q_0}+\frac{1}{q_1}$, with $0<s$, $1<p<+\infty$ and $1<p_0,p_1, q_0, q_1\leq +\infty$.
\item[2)] For $0<s, s_1, s_2<1$ with $s=s_1+s_2$ and $1<p, p_1, p_2<+\infty$ with $\frac{1}{p}=\frac{1}{p_1}+\frac{1}{p_2}$, we have
$$\|(-\Delta)^{\frac{s}{2}}(fg)-(-\Delta)^{\frac{s}{2}}(f)g-(-\Delta)^{\frac{s}{2}}(g)f\|_{L^p}\leq C\|(-\Delta)^{\frac{s_1}{2}}f\|_{L^{p1}}\|(-\Delta)^{\frac{s_2}{2}}g\|_{L^{p_2}}.$$
\end{itemize}
\end{Lemme}
See \cite{Naibo} and \cite{GrafakosOh} for a proof of these estimates.
In the case of the $L^2$-based Sobolev spaces we also have the following useful estimate: 
\begin{Lemme}[Product rule in Sobolev spaces]\label{ProductRule}
For $0\leq s<+\infty$ and $0<\delta<\frac32$,
$$\|fg\|_{\dot{H}^{s+\delta-\frac32}}\leq C\left(\|f\|_{\dot{H}^{\delta}}\|g\|_{\dot{H}^{s}}+\|g\|_{\dot{H}^{\delta}}\|f\|_{\dot{H}^{s}}\right).$$
\end{Lemme}
See \cite[Lemma 7.3]{PGLR2} for a proof of this inequality.
\mysection{Proof of Theorem \ref{Theo_ExistenceStationnaires}}\label{Secc_ProofTheo1}
We apply the Leray projector $\mathbb{P}(\vpsi)=\vpsi+\vn \frac{1}{(-\Delta)}(\vn\cdot \vpsi)$ to obtain on the one hand the following equation of the velocity (recall that $\dv(\vu)=\dv(\vf)=0$):
\begin{equation}\label{NSStationnaire23}
(-\Delta)^{\frac \alpha 2}  \vu +\mathbb{P}((\vu\cdot \vn)\vu)-\vf=0,
\end{equation}
and on the other hand, using the divergence free condition of $\vu$ and $\vf$, we have the equation for the pressure
\begin{equation}\label{NSStationnaire3}
p=\frac{1}{(-\Delta)}\left(\dv((\vu\cdot\vn )\vu)\right).
\end{equation}
We can thus focus our study on the velocity field $\vu$ and then we will deduce the properties needed for the pressure $p$. In order to solve equation (\ref{NSStationnaire23}) we will first consider a function $\phi\in \mathcal{C}^{\infty}_{0}(\R)$ such that $0\leq \theta (x)\leq 1$ with $\phi (x)=1$ if $|x|\leq 1$ and $\theta (x)=0$ if $|x|>2$, then for $R>1$ we set $\theta_{R}(x)=\theta(\frac{x}{R})$.
 With this auxiliar function and for some $0<\epsilon<1$ we study the following equation
\begin{equation}\label{EquationRegulariseeStationnaire}
-\epsilon\Delta \vu+(-\Delta)^{\frac \alpha 2}   \vu +\mathbb{P}\left(\left[(\theta_{R}\vu)\cdot \vn\right](\theta_{R}\vu)\right)-\vf=0,\qquad \dv(\vu)=0.
\end{equation}
Remark that, at least formally, if we make $\epsilon\to 0$ and $R\to +\infty$, we recover the equation (\ref{NSStationnaire23}).\\

The previous equation (\ref{EquationRegulariseeStationnaire}) can be seen as a perturbation of the  stationary Navier-Stokes system (\ref{NSStationnaire13}) and we will study the existence of solutions for this modified problem using the structure of the usual stationary Navier-Stokes. Indeed, we note  that this equation can be rewritten as 
\begin{equation}\label{ProblemePointFixe}
\vu=T_{R, \epsilon}(\vu),
\end{equation}
where
\begin{equation}\label{DefTREpsi}
T_{R, \epsilon}(\vu)=\frac{-1}{[-\epsilon\Delta+(-\Delta)^{\frac{\alpha}{2}}]}\left(\mathbb{P}\left(\left[(\theta_{R}\vu)\cdot \vn\right](\theta_{R}\vu)\right)-\vf\right).
\end{equation}
Thus, in order to obtain a solution for the problem $\vu=T_{R, \epsilon}(\vu)$ we will apply the Schaefer fixed point theorem (see \cite[Theorem 16.1]{PGLR2}):
\begin{Theoreme}[Schaefer]\label{Teo_Schaefer}
Consider the following functional space:
\begin{equation}\label{DefFuncSpace}
E=\big\{\vv:\R\longrightarrow \R: \vv\in \dot{H}^{1}(\R) \text{ and } \dv(\vv)=0\big\}.
\end{equation}
If we have the following points:
\begin{itemize}
\item[$1)$] the application $T_{R, \epsilon}$ defined in (\ref{DefTREpsi}) is continuous and compact in the space $E$,
\item[$2)$] if $\vu=\lambda T_{R}(\vu)$ for all $\lambda\in [0,1]$, then we have $\|\vu\|_{\dot{H}^{1}}\leq M$,
\end{itemize}
then the equation (\ref{ProblemePointFixe}) admits at least one solution $\vu\in E$.
\end{Theoreme}
As we can see, in order to obtain a solution of the modified problem (\ref{EquationRegulariseeStationnaire}), it is enough to verify the two points of the previous theorem. We decompose our study in some propositions and corollaries that will be helpful in the sequel. 
\begin{Proposition}\label{PropoContinuidad}
The application $T_{R, \epsilon}$ is continuous and compact in the space $E$.
\end{Proposition}
{\bf Proof.} We start writing
$$\|T_{R, \epsilon}(\vu)\|_{\dot{H}^{1}}=\left\|\frac{-\Delta}{[-\epsilon\Delta+(-\Delta)^{\frac{\alpha}{2}}]}\frac{-1}{(-\Delta)}\mathbb{P}\left(\left[(\theta_{R}\vu)\cdot \vn\right](\theta_{R}\vu)-\vf\right)\right\|_{\dot{H}^{1}}=\left\|\frac{-\Delta}{[-\epsilon\Delta+(-\Delta)^{\frac{\alpha}{2}}]}\mathcal{T}_{R}(\vu)\right\|_{\dot{H}^{1}},$$
where the operator $\mathcal{T}_{R}(\vu)$ is given by 
\begin{equation}\label{DefTR}
\mathcal{T}_{R}(\vu)=\frac{-1}{(-\Delta)}\mathbb{P}\left(\left[(\theta_{R}\vu)\cdot \vn\right](\theta_{R}\vu)-\vf\right).
\end{equation}
Observe now that the symbol $\sigma_\epsilon$ associated to the operator $\frac{-\Delta}{[-\epsilon\Delta+(-\Delta)^{\frac{\alpha}{2}}]}$ is $\sigma_\epsilon(\xi)=\frac{|\xi|^2}{\epsilon|\xi|^2+|\xi|^\alpha}$, which is a bounded Fourier multiplier \emph{i.e.} we have the uniform estimate $\sigma_{\epsilon}(\xi)\leq \frac{C}{\epsilon}$ so we can write 
$$\|T_{R, \epsilon}(\vu)\|_{\dot{H}^{1}}=\left\|\frac{-\Delta}{[-\epsilon\Delta+(-\Delta)^{\frac{\alpha}{2}}]}\mathcal{T}_{R}(\vu)\right\|_{\dot{H}^{1}}\leq \frac{C}{\epsilon}\|\mathcal{T}_{R}(\vu)\|_{\dot{H}^{1}},
$$
but from the proof of the Theorem 16.2 in \cite{PGLR2}, we know that the operator $\mathcal{T}_{R}(\vu)$ is a continuous and compact operator in the space $E$ (recall that we have the hypothesis $\vf\in \dot{H}^{-1}(\R)$) and we can deduce from this fact that the operator $T_{R, \epsilon}$ is itself continuous and compact in the space $E$.\hfill $\blacksquare$\\

We need now to establish some additional estimates.
\begin{Proposition}\label{PropoEstimaciones0}
If $\vu$ belongs to the functional space $E$ given in (\ref{DefFuncSpace}) and if $\vu$ satisfies
$$\vu=\frac{-1}{[-\epsilon\Delta+(-\Delta)^{\frac{\alpha}{2}}]}\left(\mathbb{P}\left(\left[(\theta_{R}\vu)\cdot \vn\right](\theta_{R}\vu)\right)-\vf\right),$$
then we have $\vu\in \dot{H}^1(\R)\cap \dot{H}^{\frac{\alpha}{2}}(\R)$, with $0<\alpha< 2$.
\end{Proposition}
{\bf Proof.} By the previous proposition we already know that if $\vu\in \dot{H}^1(\R)$ then the quantity $$\frac{-1}{[-\epsilon\Delta+(-\Delta)^{\frac{\alpha}{2}}]}\left(\mathbb{P}\left(\left[(\theta_{R}\vu)\cdot \vn\right](\theta_{R}\vu)\right)-\vf\right),$$ 
also belongs to $\dot{H}^1(\R)$: we only need to study if $\vu\in \dot{H}^{\frac{\alpha}{2}}(\R)$. We thus write
$$\|\vu\|_{\dot{H}^{\frac{\alpha}{2}}}=\|(-\Delta)^{\frac{\alpha}{4}}\vu\|_{L^2}=\left\|\frac{(-\Delta)^{\frac{1}{2}+\frac{\alpha}{4}}}{[-\epsilon\Delta+(-\Delta)^{\frac{\alpha}{2}}]}\frac{(-\Delta)^{\frac{1}{2}}}{(-\Delta)}\left(-\mathbb{P}\left(\left[(\theta_{R}\vu)\cdot \vn\right](\theta_{R}\vu)\right)+\vf\right)\right\|_{L^2},$$
note that the symbol $\widetilde{\sigma}_\epsilon(\xi)=\frac{|\xi|^{1+\frac{\alpha}{2}}}{\epsilon|\xi|^2+|\xi|^\alpha}$ is a bounded Fourier multiplier as we have $\widetilde{\sigma}_\epsilon(\xi)\leq \frac{C}{\epsilon}$ and we write
\begin{eqnarray*}
\|\vu\|_{\dot{H}^{\frac{\alpha}{2}}}&\leq& \frac{C}{\epsilon}\left\|\frac{(-\Delta)^{\frac{1}{2}}}{(-\Delta)}\left(-\mathbb{P}\left(\left[(\theta_{R}\vu)\cdot \vn\right](\theta_{R}\vu)\right)+\vf\right)\right\|_{L^2}=\frac{C}{\epsilon}\left\|\frac{-1}{(-\Delta)}\left(\mathbb{P}\left(\left[(\theta_{R}\vu)\cdot \vn\right](\theta_{R}\vu)\right)-\vf\right)\right\|_{\dot{H}^1}\\
&\leq & \frac{C}{\epsilon}\|\mathcal{T}_{R}(\vu)\|_{\dot{H}^{1}},
\end{eqnarray*}
where we used the definition of the operator $\mathcal{T}_{R}$ given in (\ref{DefTR}) above. We recall now that the operator $\mathcal{T}_{R}$ is bounded in the space $\dot{H}^1(\R)$ (recall that we have $\vf\in \dot{H}^{-1}(\R)$) then, as we are assuming that $\vu \in \dot{H}^1(\R)$, we have:
$$\|\vu\|_{\dot{H}^{\frac{\alpha}{2}}}\leq \frac{C_R}{\epsilon}\|\vu\|_{\dot{H}^{1}}\|\vu\|_{\dot{H}^{1}}<+\infty.$$
We have thus proven that  $\vu\in \dot{H}^1(\R)\cap \dot{H}^{\frac{\alpha}{2}}(\R)$.  \hfill $\blacksquare$\\

This proposition shows us that, although the operator $T_{R, \epsilon}$ defined in (\ref{DefTREpsi}) is bounded in the space $\dot{H}^1(\R)$, we have some additional boundedness properties in the space $\dot{H}^{\frac{\alpha}{2}}(\R)$ with $0<\alpha<2$.

\begin{Proposition}\label{PropoEstimaciones}
Let $0\leq \lambda\leq 1$. If $\vu$ belongs to the functional space $E$ given in (\ref{DefFuncSpace})  and if $\vu$ satisfies
\begin{equation}\label{EquationLambda}
\vu=\lambda\left[\frac{-1}{[-\epsilon\Delta+(-\Delta)^{\frac{\alpha}{2}}]}\left(\mathbb{P}\left(\left[(\theta_{R}\vu)\cdot \vn\right](\theta_{R}\vu)\right)-\vf\right)\right],
\end{equation}
for $0<\alpha<2$, then we have the inequality 
\begin{equation}\label{Inegalite_Energie1}
\epsilon\|\vu\|_{\dot{H}^1}^2+\|\vu\|_{\dot{H}^{\frac{\alpha}{2}}}^2\leq \lambda \|\vu\|_{\dot{H}^{\frac{\alpha}{2}}}\|\vf\|_{\dot{H}^{-\frac{\alpha}{2}}}.
\end{equation}
\end{Proposition}
{\bf Proof.} Let us first remark that since we are working in the space $E$, we have enough regularity to show that $\mathbb{P}\left(\left[(\theta_{R}\vu)\cdot \vn\right](\theta_{R}\vu)\right)\in \dot{H}^{-1}(\R)$. Indeed, we write by the properties of the Leray projector and by the Sobolev embedding $ \dot{H}^{-1}(\R)\subset L^{\frac65}(\R)$:
$$\left\|\mathbb{P}\left(\left[(\theta_{R}\vu)\cdot \vn\right](\theta_{R}\vu)\right)\right\|_{\dot{H}^{-1}}\leq C\left\|\left[(\theta_{R}\vu)\cdot \vn\right](\theta_{R}\vu)\right\|_{\dot{H}^{-1}}\leq C\left\|\left[(\theta_{R}\vu)\cdot \vn\right](\theta_{R}\vu)\right\|_{L^{\frac65}},$$
now, by the H\"older inequalities we obtain
\begin{eqnarray*}
&\leq & C\sum_{j=1}^3\left\|(\theta_{R}u_j)\partial_j(\theta_{R}\vu)\right\|_{L^{\frac65}}\leq C\sum_{j=1}^3\|\theta_{R}u_j\|_{L^3}\|\partial_j(\theta_{R}\vu)\|_{L^2}\\
&\leq & C\sum_{j=1}^3\|\theta_{R}\|_{L^6}\|u_j\|_{L^6}\left(\|(\partial_j\theta_{R})\vu\|_{L^2}+\|\theta_{R}(\partial_j\vu)\|_{L^2}\right)\\
&\leq& C\sum_{j=1}^3\|\theta_{R}\|_{L^6}\|u_j\|_{L^6}\left(\|\partial_j\theta_{R}\|_{L^3}\|\vu\|_{L^6}+\|\theta_{R}\|_{L^\infty}\|\partial_j\vu\|_{L^2}\right)\leq  C_R\|\vu\|_{L^6}(\|\vu\|_{L^6}+\|\vu\|_{\dot{H}^1})\\
&\leq &C_R\|\vu\|_{\dot{H}^1}\|\vu\|_{\dot{H}^1}<+\infty,
\end{eqnarray*}
where we used the Sobolev embedding $\dot{H}^1(\R)\subset L^{6}(\R)$ in the last estimate above. With this information at hand and since $\dv(\vu)=0$ we can write, by the properties of the Leray projector:
$$\int_{\R}\vu\cdot \mathbb{P}\left(\left[(\theta_{R}\vu)\cdot \vn\right](\theta_{R}\vu)\right)dx=\int_{\R}\vu\cdot \left(\left[(\theta_{R}\vu)\cdot \vn\right](\theta_{R}\vu)\right)dx,$$
but since by an integration by parts we have 
$$\int_{\R}\vu\cdot \left(\left[(\theta_{R}\vu)\cdot \vn\right](\theta_{R}\vu)\right)dx=-\int_{\R}\vu\cdot \left(\left[(\theta_{R}\vu)\cdot \vn\right](\theta_{R}\vu)\right)dx,$$
we deduce that 
\begin{equation}\label{IdentiteNulle}
\int_{\R}\vu\cdot \left(\left[(\theta_{R}\vu)\cdot \vn\right](\theta_{R}\vu)\right)dx=0.
\end{equation}
With this information, we rewrite now the equation (\ref{EquationLambda}) in the following form
$$[-\epsilon\Delta+(-\Delta)^{\frac{\alpha}{2}}]\vu=-\lambda\left[\left(\mathbb{P}\left(\left[(\theta_{R}\vu)\cdot \vn\right](\theta_{R}\vu)\right)-\vf\right)\right]=-\lambda\mathbb{P}\left(\left[(\theta_{R}\vu)\cdot \vn\right](\theta_{R}\vu)\right)+\lambda\mathbb{P}\left(\vf\right),$$
from which we deduce 
$$-\epsilon\int_{\R}\Delta\vu \cdot \vu\, dx+\int_{\R}(-\Delta)^{\frac{\alpha}{2}}\vu\cdot \vu\,dx=-\lambda\int_{\R}\mathbb{P}\left(\left[(\theta_{R}\vu)\cdot \vn\right](\theta_{R}\vu)\right)\cdot \vu\,dx+\lambda\int_{\R}\mathbb{P}\left(\vf\right)\cdot \vu\,dx.$$
Using the identity (\ref{IdentiteNulle}), by the properties of the Leray projector, since $\dv(\vu)=0$, and using the properties of the operators $\Delta$ and $(-\Delta)^{\frac{\alpha}{2}}$, we obtain
$$\epsilon\|\vu\|_{\dot{H}^1}^2+\|\vu\|_{\dot{H}^{\frac{\alpha}{2}}}^2=\lambda\int_{\R}(-\Delta)^{-\frac{\alpha}{4}}\vf\cdot(-\Delta)^{\frac{\alpha}{4}} \vu\,dx\leq \lambda\|\vf\|_{\dot{H}^{-\frac{\alpha}{2}}}\|\vu\|_{\dot{H}^{\frac{\alpha}{2}}},$$
and we have proven the estimate (\ref{Inegalite_Energie1}). \hfill $\blacksquare$\\

This estimate has several consequences and we gather them in the following corollary:
\begin{Corollaire}\label{CoroEnergieEstimate}
In the general framework of the Proposition \ref{PropoEstimaciones}, \emph{i.e.} if $\vu$ belongs to the functional space $E$ given in (\ref{DefFuncSpace}) and if $\vu$ satisfies the equation (\ref{EquationLambda}), then we have the following points:
\begin{itemize}
\item[1)] For $\epsilon>0$ one has the inequality
\begin{equation}\label{Inegalite_Energie2}
\|\vu\|_{\dot{H}^1}\leq \frac{1}{\sqrt{2\epsilon}}\|\vf\|_{\dot{H}^{-\frac{\alpha}{2}}}.
\end{equation}
\item[2)] We also have the uniform estimate 
\begin{equation}\label{Inegalite_Energie3}
\|\vu\|_{\dot{H}^{\frac{\alpha}{2}}}\leq \|\vf\|_{\dot{H}^{-\frac{\alpha}{2}}}.
\end{equation}
\end{itemize}
\end{Corollaire}
{\bf Proof.} From the estimate (\ref{Inegalite_Energie1}) we write, by the Young inequalities for the product:
$$\epsilon\|\vu\|_{\dot{H}^1}^2+\|\vu\|_{\dot{H}^{\frac{\alpha}{2}}}^2\leq \frac{\lambda}{2}\left(\|\vf\|_{\dot{H}^{-\frac{\alpha}{2}}}^2+\|\vu\|_{\dot{H}^{\frac{\alpha}{2}}}^2\right),$$
from which we easily obtain $\epsilon\|\vu\|_{\dot{H}^1}^2+(1-\tfrac{\lambda}{2})\|\vu\|_{\dot{H}^{\frac{\alpha}{2}}}^2\leq \frac{\lambda}{2}\|\vf\|_{\dot{H}^{-\frac{\alpha}{2}}}^2$ and since $0\leq \lambda\leq 1$ we easily deduce the two wished estimates.  \hfill $\blacksquare$\\

\noindent {\bf End of the proof of the Theorem \ref{Theo_ExistenceStationnaires}.} By the Schaefer fixed-point theorem, in order to obtain the existence of a solution of the problem (\ref{ProblemePointFixe})-(\ref{DefTREpsi}) we only need to prove the two points given in Theorem \ref{Teo_Schaefer}. Thus, for some fixed parameters $R>1$, $\epsilon>0$ we know by the Proposition \ref{PropoContinuidad} that the application $T_{R, \epsilon}$ is continuous and compact in the space $E$ given in (\ref{DefFuncSpace}). The second point of Theorem \ref{Teo_Schaefer} is given by the estimate (\ref{Inegalite_Energie2}) stated in the Corollary \ref{CoroEnergieEstimate}. We have thus the existence of a solution of the problem
$$\vu=\frac{-1}{[-\epsilon\Delta+(-\Delta)^{\frac{\alpha}{2}}]}\left(\mathbb{P}\left(\left[(\theta_{R}\vu)\cdot \vn\right](\theta_{R}\vu)\right)-\vf\right),$$ 
where $\vu\in \dot{H}^1(\R)\cap \dot{H}^{\frac{\alpha}{2}}(\R)$. 
\begin{Remarque}
Note that the solution obtained above depends on the parameters $R>1$ and $\epsilon>0$ and they will be denoted by $\vu_{R,\epsilon}$.
\end{Remarque}
Now we need to recover the initial problem by making $R\to +\infty$ and $\epsilon\to 0$ in the solutions $\vu_{R, \epsilon}$. To do so, we will first fix $\epsilon>0$ and then we will take the limit when  $R\to +\infty$. Indeed, we remark that for a fixed $\epsilon>0$, we have the uniform (in $R$) estimate
(\ref{Inegalite_Energie3}), and thus there exists a sequence $R_k\to+\infty$ such that $\vu_{R_k, \epsilon}$ converges weakly in $\dot{H}^{\frac{\alpha}{2}}(\R)$ to some limit $\vu_{\epsilon}$. Moreover, by the Lemma \ref{Lemma_Rellich-Kondrachov}, we have for all $0<\alpha<2$ the strong convergence of $\vu_{R_k, \epsilon}$ to some limit $\vu_{\epsilon}$ in the space $L^2_{loc}(\R)$. These two facts allows us to obtain a weak convergence (in $\mathcal{D}'$) of the nonlinear term $\left(\left[(\theta_{R}\vu_{R, \epsilon})\cdot \vn\right](\theta_{R}\vu_{R, \epsilon})\right)$ to $\left(\left[\vu_{\epsilon}\cdot \vn\right]\vu_{\epsilon}\right)$ when $R\to +\infty$. We thus obtain a function $\vu_{\epsilon}$ which is a solution of the problem
$$-\epsilon\Delta \vu_{\epsilon}+(-\Delta)^{\frac \alpha 2}   \vu_{\epsilon} +\mathbb{P}\left(\left[\vu_{\epsilon}\cdot \vn\right]\vu_{\epsilon}\right)-\vf=0.$$
Similarly, since we have the uniform (in $\epsilon$) control $\|\vu_{\epsilon}\|_{\dot{H}^{\frac{\alpha}{2}}}\leq \|\vf\|_{\dot{H}^{-\frac{\alpha}{2}}}$ given in the estimate (\ref{Inegalite_Energie3}), there exists a subsequence $\epsilon_k\to 0$ such that $\vu_{\epsilon_k}$ convergence weakly to a limit $\vu$ in the space $\dot{H}^{\frac{\alpha}{2}}(\R)$. Again, by the Lemma \ref{Lemma_Rellich-Kondrachov}, we obtain the strong convergence of $\vu_{\epsilon_k}$ to $\vu$ in $L^2_{loc}(\R)$ and from these facts we obtain the weak convergence (in $\mathcal{D}'$) of the quantity $\left[\vu_{\epsilon}\cdot \vn\right]\vu_{\epsilon}$ to $(\vu\cdot \vn)\vu$ when $\epsilon\to 0$. We have thus obtained a solution $\vu$ of the equation
$$(-\Delta)^{\frac \alpha 2} \vu+\mathbb{P}\left((\vu\cdot \vn)\vu\right)-\vf=0,$$
which belongs to the space $\dot{H}^{\frac{\alpha}{2}}(\R)$ and satisfies $\|\vu\|_{\dot{H}^{\frac{\alpha}{2}}}\leq \|\vf\|_{\dot{H}^{-\frac{\alpha}{2}}}$.\\

To end the proof, we need to study the pressure $p$. By the equation (\ref{NSStationnaire3}) and by the divergence free property of $\vu$ we can write
$$\|p\|_{\dot{H}^{\alpha-\frac32}}=\left\|(-\Delta)^{\frac{\alpha-\frac{3}{2}}{2}}p\right\|_{L^2}=\left\|\frac{(-\Delta)^{\frac{\alpha-\frac{3}{2}}{2}}}{(-\Delta)}\dv(\dv(\vu\otimes\vu))\right\|_{L^2}=\left\|(-\Delta)^{\frac{\alpha-\frac{3}{2}}{2}}(\vu\otimes\vu)\right\|_{L^2}=\|\vu\otimes \vu\|_{\dot{H}^{\alpha-\frac32}}.$$
Now by the Product rule in Sobolev spaces given in the Lemma \ref{ProductRule} we have
$$\|\vu\otimes \vu\|_{\dot{H}^{\alpha-\frac32}}\leq C\|\vu\|_{\dot{H}^{\frac{\alpha}{2}}}\|\vu\|_{\dot{H}^{\frac{\alpha}{2}}}<+\infty,$$
from which we easily deduce that $\|p\|_{\dot{H}^{\alpha-\frac32}}<+\infty$ and this ends the proof of the Theorem \ref{Theo_ExistenceStationnaires}.\hfill $\blacksquare$
\mysection{Proof of Theorem \ref{Theo_Liouville}}\label{Secc_ProofTheo2}
We start the proof of this theorem with a brief remark:
\begin{Lemme}\label{Lemme_InformationPression}
Let $(\vu, p)$ be a solution of the fractional Navier-Stokes equation
$$(-\Delta)^{\frac \alpha 2} \vu +(\vu\cdot \vn)\vu+\vn p=0,\qquad \dv(\vu)=0.$$
If we have $\vu\in L^q(\R)$ for some $2<q<+\infty$ then the pressure $p$ belongs to the space $L^{\frac{q}{2}}(\R)$.
\end{Lemme}
{\bf Proof.} Applying the divergence operator to the equation above and using the divergence free property of $\vu$ we obtain $\dv((\vu\cdot \vn)\vu)+\dv(\vn p)=0$ which leads us to the equation
$\Delta p=-\dv(\dv(\vu\otimes \vu))$ from which we deduce the expression $p=\frac{1}{(-\Delta)}\dv(\dv(\vu\otimes \vu))$. Thus, taking the $L^{\frac{q}{2}}(\R)$ norm, we have
$$\|p\|_{L^{\frac{q}{2}}}=\left\|\frac{1}{(-\Delta)}\dv(\dv(\vu\otimes \vu))\right\|_{L^{\frac{q}{2}}}.$$
Since $1<\frac{q}{2}<+\infty$, the Riesz transforms are bounded in the space $L^{\frac{q}{2}}(\R)$ and we can write
$$\left\|\frac{1}{(-\Delta)}\dv(\dv(\vu\otimes \vu))\right\|_{L^{\frac{q}{2}}}\leq C\|\vu\otimes \vu\|_{L^{\frac{q}{2}}}\leq C\| \vu\|_{L^{q}}\| \vu\|_{L^{q}}<+\infty,$$
and we obtain that $p\in L^{\frac{q}{2}}(\R)$. \hfill $\blacksquare$\\

This simple remark allows us to deduce some integrability results for the pressure from the information available on the velocity field $\vu$.\\

We will prove now that in the framework of the Theorem \ref{Theo_Liouville}, the unique solution of the equation (\ref{NSStationnaire2}) is the trivial solution. For this we consider $\theta \in \mathcal{C}^{\infty}_0(\R)$ a smooth cut-off function given by $0\leq \theta \leq 1$, $\theta(x)=1$ if $\vert x \vert <\frac{1}{2}$ and $\theta(x)=0$ if $\vert x \vert \geq 1$. For $R>1$ a real parameter, we define the function 
$$\theta_R(x)=\theta\left( \frac{x}{R}\right),$$ 
in particular we have $\theta_R(x)=1$ if $\vert x \vert < \frac{R}{2}$ and $\theta_R(x)=0$ if $\vert x \vert \geq R$ and thus $supp(\theta_R)\subset B_R$, where $B_R$ denotes the ball $B(0,R)$. With this auxiliary function, we multiply the equation (\ref{NSStationnaire2}) by $(\theta_R \vu)$ and we integrate:
\begin{equation}\label{Identite}
\int_{\R} \underbrace{(-\Delta)^{\frac{\alpha}{2}} \vu\cdot \left( \theta_R \vu\right)}_{(1)} +\underbrace{(\vu \cdot \vec{\nabla})\vu\cdot \left( \theta_R \vu\right)}_{(2)}  + \underbrace{\vec{\nabla}p \cdot \left( \theta_R \vu\right)}_{(3)} dx=0,
\end{equation} 
and we study each one of these terms separately. For the first term in (\ref{Identite}) we write, using the properties of the operator $(-\Delta)^{\frac{\alpha}{2}}$:
\begin{eqnarray*}
\int_{\R} (-\Delta)^{\frac{\alpha}{2}} \vu\cdot \left( \theta_R \vu\right)dx&=&\int_{\R} (-\Delta)^{\frac{\alpha}{4}} \vu\cdot (-\Delta)^{\frac{\alpha}{4}}\left( \theta_R \vu\right)dx\\
&=&\int_{\R} (-\Delta)^{\frac{\alpha}{4}} \vu\cdot \left[((-\Delta)^{\frac{\alpha}{4}}\vu) \theta_R + (-\Delta)^{\frac{\alpha}{4}}( \theta_R \vu)-((-\Delta)^{\frac{\alpha}{4}}\vu) \theta_R\right]dx,
\end{eqnarray*}
and we have
\begin{equation}\label{Identite1}
\int_{\R} (-\Delta)^{\frac{\alpha}{2}} \vu\cdot \left( \theta_R \vu\right)dx
=\int_{B_R} |(-\Delta)^{\frac{\alpha}{4}} \vu|^2\theta_R dx+ \int_{\R} (-\Delta)^{\frac{\alpha}{4}} \vu\cdot\left[(-\Delta)^{\frac{\alpha}{4}}( \theta_R \vu)-((-\Delta)^{\frac{\alpha}{4}}\vu) \theta_R\right]dx,
\end{equation}
where we used the fact that $supp(\theta_R)\subset B_R$ in the second integral above. \\[5mm]
For the second term of (\ref{Identite}) we have:
\begin{eqnarray*} 
\int_{\R} (\vu \cdot \vec{\nabla})\vu\cdot ( \theta_R \vu )dx &=& \sum_{i,j=1}^{3} \int_{\R} u_j (\partial_{x_j} u_i) ( \theta_R u_i ) dx =  \sum_{i,j=1}^{3} \int_{\R} \theta_R u_j  (\partial_{x_j} u_i) u_i dx\\
&=&  \sum_{i,j=1}^{3} \int_{\R} \theta_R u_j  (\partial_{x_j} \left(\frac{u^{2}_{i}}{2}\right)) dx,
\end{eqnarray*}
and by an integration by parts we obtain
$$\sum_{i,j=1}^{3} \int_{\R} \theta_R u_j  (\partial_{x_j} \left(\frac{u^{2}_{i}}{2}\right)) dx =-\sum_{i,j=1}^{3} \int_{\R} \theta_R \left(\partial_{x_j}  u_j\right) \frac{u^{2}_{i}}{2} dx  - \int_{\R} \vec{\nabla} \theta_R \cdot \left( \frac{\vert \vu\vert^2}{2} \vu\right)dx.$$
Now, using the fact that $\dv(\vu)=0$, we have that the second integral above is null and we can write
\begin{equation} \label{Identite2}
\int_{\R} (\vu \cdot \vec{\nabla})\vu\cdot ( \theta_R \vu )dx=- \int_{B_R} \vec{\nabla} \theta_R \cdot \left( \frac{\vert \vu\vert^2}{2} \vu\right)dx,
\end{equation}
where we used the support property of the auxiliar function $\theta_R$. \\[5mm]
Finally, for the last term of (\ref{Identite}), by an integration by parts, using again the fact $\dv(\vu)=0$ and the support property of $\theta_R$, we obtain
\begin{eqnarray} \nonumber
\int_{\R} \vec{\nabla}p\cdot ( \theta_R \vu )dx&=& 	\sum_{i=1}^{3}\int_{\R}(\partial_{x_i} p) \theta_R u_i dx=- \sum_{i=1}^{3}\int_{\R} p \,\partial_{x_i} (\theta_R u_i)dx \\
&=& - \sum_{i=1}^{3}\int_{\R} p (\partial_{x_i} \theta_R)  ( u_i) dx = - \int_{B_R}\vec{\nabla}\theta_R \cdot(p \vu) dx. \label{Identite3}
\end{eqnarray}
Thus, with the expressions  (\ref{Identite1}), (\ref{Identite2}), and (\ref{Identite3}), we can rewrite equation (\ref{Identite}) in the following manner:
\begin{eqnarray*}
\int_{B_R} |(-\Delta)^{\frac{\alpha}{4}} \vu|^2\theta_R dx&+& \int_{\R} (-\Delta)^{\frac{\alpha}{4}} \vu\cdot\left[(-\Delta)^{\frac{\alpha}{4}}( \theta_R \vu)-((-\Delta)^{\frac{\alpha}{4}}\vu) \theta_R\right]dx\\
&&-\int_{B_R} \vec{\nabla} \theta_R \cdot \left( \frac{\vert \vu\vert^2}{2} \vu\right)dx - \int_{B_R}\vec{\nabla}\theta_R \cdot(p \vu) dx=0,
\end{eqnarray*}
from which we obtain the equation
\begin{eqnarray*}
\int_{B_R} |(-\Delta)^{\frac{\alpha}{4}} \vu|^2\theta_R dx&=& \int_{\R} (-\Delta)^{\frac{\alpha}{4}} \vu\cdot\left[((-\Delta)^{\frac{\alpha}{4}}\vu) \theta_R-(-\Delta)^{\frac{\alpha}{4}}( \theta_R \vu)\right]dx+\int_{B_R} \vec{\nabla} \theta_R \cdot \left( \frac{\vert \vu\vert^2}{2} \vu\right)dx\\
&& + \int_{B_R}\vec{\nabla}\theta_R \cdot(p \vu) dx.
\end{eqnarray*}
We recall now that since $0\leq \theta_R(x)\leq 1$ and $\theta_R (x)=1$ if $\vert x \vert < \frac{R}{2}$,  we have the estimate
$$\int_{B_{\frac{R}{2}}} |(-\Delta)^{\frac{\alpha}{4}} \vu|^2 dx \leq  \int_{B_R}|(-\Delta)^{\frac{\alpha}{4}} \vu|^2\theta_R dx,$$ 
and we can write
\begin{eqnarray}
\int_{B_{\frac{R}{2}}} |(-\Delta)^{\frac{\alpha}{4}} \vu|^2 dx &\leq & \underbrace{\int_{\R} (-\Delta)^{\frac{\alpha}{4}} \vu\cdot\left[((-\Delta)^{\frac{\alpha}{4}}\vu) \theta_R-(-\Delta)^{\frac{\alpha}{4}}( \theta_R \vu)\right]dx}_{(I_a)}+\underbrace{\int_{B_R} \vec{\nabla} \theta_R \cdot \left( \frac{\vert \vu\vert^2}{2} \vu\right)dx}_{(I_b)}\notag\\
&& + \underbrace{\int_{B_R}\vec{\nabla}\theta_R \cdot(p \vu)dx}_{(I_c)}.\label{Inegalite}
\end{eqnarray}
We will now show that we have $\underset{R\to+\infty}{\lim}I_a=0$, $\underset{R\to+\infty}{\lim}I_b=0$ and $\underset{R\to+\infty}{\lim}I_c=0$. Indeed:

\begin{itemize}

\item[(a)] For the first term above, we write by the Cauchy-Schwarz inequality
\begin{eqnarray*}
I_a&=&\int_{\R}(-\Delta)^{\frac{\alpha}{4}} \vu\cdot\left[((-\Delta)^{\frac{\alpha}{4}}\vu) \theta_R-(-\Delta)^{\frac{\alpha}{4}}( \theta_R \vu)\right]dx\\
&\leq&\|(-\Delta)^{\frac{\alpha}{4}} \vu\|_{L^2}\left\|
((-\Delta)^{\frac{\alpha}{4}}\vu) \theta_R-(-\Delta)^{\frac{\alpha}{4}}( \theta_R \vu)\right\|_{L^2}\\
&\leq &\|\vu\|_{\dot{H}^{\frac{\alpha}{2}}}\left(\left\|(-\Delta)^{\frac{\alpha}{4}}\big( \theta_R \vu\big)-
\big((-\Delta)^{\frac{\alpha}{4}}\vu\big) \theta_R-\big((-\Delta)^{\frac{\alpha}{4}} \theta_R\big) \vu\right\|_{L^2}+\left\|\big((-\Delta)^{\frac{\alpha}{4}} \theta_R\big) \vu\right\|_{L^2}\right).
\end{eqnarray*}
We apply now the second point of Lemma \ref{FracLeibniz} to obtain the estimate
\begin{equation*}
I_a\leq \|\vu\|_{\dot{H}^{\frac{\alpha}{2}}}\left(\|(-\Delta)^{\frac{\alpha_1}{4}}\theta_R\|_{L^{p_1}}\|(-\Delta)^{\frac{\alpha_2}{4}}  \vu\|_{L^{p_2}}+\left\|\big((-\Delta)^{\frac{\alpha}{4}} \theta_R\big) \vu\right\|_{L^2}\right),
\end{equation*}
where $\alpha=\alpha_1+\alpha_2$, $0<\alpha, \alpha_1, \alpha_2<2$ and $\frac{1}{2}=\frac{1}{p_1}+\frac{1}{p_2}$. \\

Recall now that we have $0<\alpha<2$ and that we are assuming in all the cases stated in Theorem \ref{Theo_Liouville} the condition $\vu\in L^{\frac{6-\epsilon}{3-\alpha}}(\R)$ with $0<\epsilon<2\alpha$, thus by the H\"older inequality with $\frac{1}{2}=\frac{2\alpha-\epsilon}{12-2\epsilon}+\frac{3-\alpha}{6-\epsilon}$ we have
\begin{eqnarray*}
I_a&\leq &\|\vu\|_{\dot{H}^{\frac{\alpha}{2}}}\left(\|(-\Delta)^{\frac{\alpha_1}{4}}\theta_R\|_{L^{p_1}}\|(-\Delta)^{\frac{\alpha_2}{4}}  \vu\|_{L^{p_2}}+\|(-\Delta)^{\frac{\alpha}{4}} \theta_R\|_{L^{\frac{12-2\epsilon}{2\alpha-\epsilon}}}\|\vu\|_{L^\frac{6-\epsilon}{3-\alpha}}\right)\\
&\leq& \|\vu\|_{\dot{H}^{\frac{\alpha}{2}}}\left(CR^{-\frac{\alpha_1}{2}+\frac{3}{p_1}}\|(-\Delta)^{\frac{\alpha_2}{4}}  \vu\|_{L^{p_2}}+CR^{-\frac{\alpha}{2}+3\frac{2\alpha-\epsilon}{12-2\epsilon}}\|\vu\|_{L^\frac{6-\epsilon}{3-\alpha}}\right),
\end{eqnarray*}
where we used the properties of the function $\theta_R$ in the last estimate above. Let us also note that, due to the complex interpolation theory (see \cite[Theorem 6.4.5.]{BL}) we have 
$$\big[\dot{H}^{\frac{\alpha}{2}}, L^{\frac{6-\epsilon}{3-\alpha}}\big]_{\nu}=\dot{W}^{\frac{\alpha_2}{2}, p_2}\qquad \mbox{and}\qquad \|(-\Delta)^{\frac{\alpha_2}{4}}  \vu\|_{L^{p_2}}=\|\vu\|_{\dot{W}^{\frac{\alpha_2}{2},p_2}} \leq C\|\vu\|_{\dot{H}^{\frac{\alpha}{2}}}^\nu\|\vu\|_{L^\frac{6-\epsilon}{3-\alpha}}^{1-\nu},$$
with the relationships 
\begin{equation}\label{ConditionInterpol}
\alpha_2=\nu\alpha, \qquad \frac{1}{p_2}=\frac{\nu}{2}+(1-\nu)\frac{3-\alpha}{6-\epsilon}\qquad \mbox{for some} \qquad 0<\nu<1,
\end{equation}
and then we can write
$$I_a\leq \|\vu\|_{\dot{H}^{\frac{\alpha}{2}}}\left(CR^{-\frac{\alpha_1}{2}+\frac{3}{p_1}}\|\vu\|_{\dot{H}^{\frac{\alpha}{2}}}^\nu\|\vu\|_{L^\frac{6-\epsilon}{3-\alpha}}^{1-\nu}+CR^{-\frac{\alpha}{2}+3\frac{2\alpha-\epsilon}{12-2\epsilon}}\|\vu\|_{L^\frac{6-\epsilon}{3-\alpha}}\right).
$$
But since we have $\alpha=\alpha_1+\alpha_2$ and $\frac{1}{2}=\frac{1}{p_1}+\frac{1}{p_2}$, following the conditions (\ref{ConditionInterpol}) above, we obtain that $\alpha_1=(1-\nu)\alpha$ and $\frac{1}{p_1}=(1-\nu)\frac{2\alpha-\epsilon}{12-2\epsilon}$ and we can write 
\begin{eqnarray*}
I_a&\leq &C\|\vu\|_{\dot{H}^{\frac{\alpha}{2}}}\Big(R^{(1-\nu)[-\frac{\alpha}{2}+\frac{6\alpha-3\epsilon}{12-2\epsilon}]}\|\vu\|_{\dot{H}^{\frac{\alpha}{2}}}^\nu\|\vu\|_{L^\frac{6-\epsilon}{3-\alpha}}^{1-\nu}+R^{-\frac{\alpha}{2}+\frac{6\alpha-3\epsilon}{12-2\epsilon}}\|\vu\|_{L^\frac{6-\epsilon}{3-\alpha}}\Big).
\end{eqnarray*}
Remark that since $0< \alpha<2$ and $0<\epsilon<2\alpha$ we always have $\frac{6\alpha-3\epsilon}{12-2\epsilon}<\frac{\alpha}{2}$ and thus all the powers of the parameter $R$ in the right-hand side above are  negative. Moreover we have $\|\vu\|_{\dot{H}^{\frac{\alpha}{2}}}<+\infty$ and $\|\vu\|_{L^\frac{6-\epsilon}{3-\alpha}}<+\infty$, so we obtain 
\begin{equation}\label{EstimateR1}
\underset{R\to+\infty}{\lim}I_a=0.
\end{equation}
\begin{Remarque}
Note that in all the cases $\alpha=1$, $1<\alpha<2$ and $\frac35<\alpha<1$ stated in Theorem \ref{Theo_Liouville}, in order to obtain the previous limit (\ref{EstimateR1}) we only require the information $\vu\in L^\frac{6-\epsilon}{3-\alpha}(\R)$ for some $0<\epsilon<2\alpha$ and no further conditions are needed for the parameter $\epsilon$. The conditions (\ref{ConditionIntegrable}) and (\ref{ConditionIntegrable1}) will appear in the study of the limits $\underset{R\to+\infty}{\lim}I_b$ and $\underset{R\to+\infty}{\lim}I_c$.
\end{Remarque}

\item[(b)] For the second term of (\ref{Inegalite}) we recall that $\theta_R (x)=1$ if $\vert x \vert <\frac{R}{2}$ and $\theta_R(x)=0$ if $\vert x \vert \geq R$ and thus we have
\begin{equation}\label{ProprieteSupport}
supp \left( \vec{\nabla}\theta_R \right) \subset \big\{ x\in \R: \tfrac{R}{2} < \vert x \vert <R \big\}=\mathcal{C}(\tfrac R2, R),
\end{equation}
and with this remark we can write 
$$I_b=\int_{B_R} \vec{\nabla} \theta_R \cdot \left(\frac{|\vu|^2}{2}\vu \right) dx=\int_{\mathcal{C}(\frac R2, R)} \vec{\nabla} \theta_R \cdot \left( \frac{|\vu|^2}{2}\vu \right) dx.$$
In order to study the limit when $R\to +\infty$, we decompose our study following the values of $\alpha$ and the information available. Indeed:
\begin{itemize}
\item[$\bullet$] If $\alpha=1$, we have $\vu\in \dot{H}^{\frac12}(\R)$ and thus, by the Sobolev embeddings we also have $\vu\in L^3(\R)$, so we can write:
\begin{equation}\label{EstimateR2}
I_b\leq C \|\vn\theta_R \|_{L^{\infty}(\mathcal{C}(\frac R2, R))}\|\vu\|_{L^3(\mathcal{C}(\frac R2, R))}^3\leq  CR^{-1} \|\vu\|_{L^{3}(\mathcal{C}(\frac R2, R))}^3,
\end{equation}
from which we easily deduce that 
$$\underset{R\to+\infty}{\lim}I_b=0.$$
\item[$\bullet$] If $1<\alpha<2$, we know by hypothesis that $\vu\in L^{\frac{6-\epsilon}{3-\alpha}}(\R)$ 
and recall that we have in this case the condition (\ref{ConditionIntegrable}), \emph{i.e.} $1+\frac{\epsilon}{3}\leq \alpha\leq \frac{5}{3}+\frac{2}{9}\epsilon$. Thus, if $1+\frac{\epsilon}{3}< \alpha\leq \frac{5}{3}+\frac{2}{9}\epsilon$, by the H\"older inequality with $\frac{3\alpha-3-\epsilon}{6-\epsilon}+3(\frac{3-\alpha}{6-\epsilon})=1$, we can write
\begin{eqnarray*}
I_b&\leq &C \|\vn\theta_R \|_{L^{\frac{6-\epsilon}{3\alpha-3-\epsilon}}(\mathcal{C}(\frac R2, R))}\|\vu\|_{L^{\frac{6-\epsilon}{3-\alpha}}(\mathcal{C}(\frac R2, R))}^3\notag\\
&\leq & CR^{-1+3\frac{3\alpha-3-\epsilon}{6-\epsilon}} \|\vu\|_{L^{\frac{6-\epsilon}{3-\alpha}}(\mathcal{C}(\frac R2, R))}^3.
\end{eqnarray*}
Then if $1+\frac{\epsilon}{3}< \alpha< \frac{5}{3}+\frac{2}{9}\epsilon$ the power of the parameter $R$ above is negative and then the quantity above will tend to $0$ if $R\to+\infty$. But if $\alpha=\frac{5}{3}+\frac{2}{9}\epsilon$, we have $-1+3\frac{3\alpha-3-\epsilon}{6-\epsilon}=0$, then since $\vu\in L^{\frac{6-\epsilon}{3-\alpha}}(\R)$, we will have $\|\vu\|_{L^{\frac{6-\epsilon}{3-\alpha}}(\mathcal{C}(\frac R2, R))}\underset{R\to +\infty}{\longrightarrow 0}$. Finally if $1+\frac{\epsilon}{3}=\alpha$ then we have $L^{\frac{6-\epsilon}{3-\alpha}}(\R)=L^{3}(\R)<+\infty$, moreover the power of $R$ is negative  and equal to $-1$ and we can proceed as in (\ref{EstimateR2}). Thus, in any case we obtain
$$\underset{R\to+\infty}{\lim}I_b=0.$$
Note that in the case $1<\alpha<2$ besides the condition $0<\epsilon<2\alpha$ we need the relationship (\ref{ConditionIntegrable}) between $\alpha$ and $\epsilon$.\\
\item[$\bullet$] If $\frac{3}{5}<\alpha<1$, in this case we have the additional condition $\vu\in L^{\frac{6+\epsilon}{3-\alpha}}(\R)$ with the relationship $1-\frac{\epsilon}{3}\leq \alpha\leq \frac{5}{3}-\frac{2}{9}\epsilon$ (recall the condition (\ref{ConditionIntegrable1})). As above, if $1-\frac{\epsilon}{3}< \alpha\leq \frac{5}{3}-\frac{2}{9}\epsilon$, by the H\"older inequality with $\frac{3\alpha-3+\epsilon}{6+\epsilon}+3(\frac{3-\alpha}{6+\epsilon})=1$, we obtain
$$I_b\leq CR^{-1+3\frac{3\alpha-3+\epsilon}{6+\epsilon}} \|\vu\|_{L^{\frac{6+\epsilon}{3-\alpha}}(\mathcal{C}(\frac R2, R))}^3.$$
Note that if $1-\frac{\epsilon}{3}< \alpha< \frac{5}{3}-\frac{2}{9}\epsilon$, the power of the parameter $R$ is negative, while if $\alpha=\frac{5}{3}-\frac{2}{9}\epsilon$ we have $-1+3\frac{3\alpha-3+\epsilon}{6+\epsilon}=0$, but we have $\|\vu\|_{L^{\frac{6+\epsilon}{3-\alpha}}(\mathcal{C}(\frac R2, R))}\underset{R\to +\infty}{\longrightarrow 0}$. Remark also that if $\alpha=1-\frac{\epsilon}{3}$, then $L^{\frac{6+\epsilon}{3-\alpha}}(\R)=L^{3}(\R)<+\infty$, the power of $R$ is equal $-1$ and we can  proceed as in (\ref{EstimateR2}). In any case we have
$$\underset{R\to+\infty}{\lim}I_b=0.$$
\begin{Remarque}
Note that when $\frac{3}{5}<\alpha<1$ we need the information $\vu\in L^{\frac{6-\epsilon}{3-\alpha}}(\R)$ with the condition $0<\epsilon<2\alpha$ in order to obtain the limit (\ref{EstimateR1}) for the term $I_a$, but we also need the information $\vu\in L^{\frac{6+\epsilon}{3-\alpha}}(\R)$ with the constraint (\ref{ConditionIntegrable1}) to obtain that $\underset{R\to+\infty}{\lim}I_b=0$.\\
Note also that the lower limit $\frac{3}{5}<\alpha$ is a consequence of the conditions $1-\frac{\epsilon}{3}\leq\alpha$ and $0<\epsilon<2\alpha$. We recall that these conditions are technical and we do not claim any optimality on them. 
\end{Remarque}
\end{itemize}
\item[(c)] For the last term of (\ref{Inegalite}) we write, using the support property (\ref{ProprieteSupport}):
$$I_c=\int_{B_R}\vec{\nabla}\theta_R \cdot(p \vu) dx=\int_{\mathcal{C}(\frac R2, R)}\vec{\nabla}\theta_R \cdot(p \vu) dx.$$
The study of this term is very similar of the previous one since by Lemma \ref{Lemme_InformationPression} we also have some information over the pressure $p$. Indeed:
\begin{itemize}
\item[$\bullet$] If $\alpha=1$, we have $\vu\in \dot{H}^{\frac12}(\R)$ and thus, by the Sobolev embeddings we have $\vu\in L^3(\R)$ but we also have $p\in L^{\frac{3}{2}}(\R)$ by Lemma \ref{Lemme_InformationPression}, and we write
\begin{equation}\label{EstimateR3}
I_c\leq C \|\vn\theta_R \|_{L^{\infty}(\mathcal{C}(\frac R2, R))}\|p\|_{L^{\frac{3}{2}}(\mathcal{C}(\frac R2, R))}\|\vu\|_{L^3(\mathcal{C}(\frac R2, R))}\leq  CR^{-1} \|p\|_{L^{\frac{3}{2}}(\R)}\|\vu\|_{L^{3}(\R)},
\end{equation}
from which we easily deduce that 
$$\underset{R\to+\infty}{\lim}I_c=0.$$
\item[$\bullet$] If $1<\alpha<2$, we have $\vu\in L^{\frac{6-\epsilon}{3-\alpha}}(\R)$ and by Lemma \ref{Lemme_InformationPression} we have $p\in L^{\frac{6-\epsilon}{6-2\alpha}}(\R)$. If $1+\frac{\epsilon}{3}< \alpha\leq \frac{5}{3}+\frac{2}{9}\epsilon$, by the H\"older inequality with $\frac{3\alpha-3-\epsilon}{6-\epsilon}+\frac{6-2\alpha}{6-\epsilon}+\frac{3-\alpha}{6-\epsilon}=1$ we obtain
\begin{eqnarray*}
I_c&\leq &\|\vn\theta_R \|_{L^{\frac{6-\epsilon}{3\alpha-3-\epsilon}}(\mathcal{C}(\frac R2, R))}\|p\|_{L^{\frac{6-\epsilon}{6-2\alpha}}(\mathcal{C}(\frac R2, R))}\|\vu\|_{L^{\frac{6-\epsilon}{3-\alpha}}(\mathcal{C}(\frac R2, R))}\notag\\
&\leq &C R^{-1+3\frac{3\alpha-3-\epsilon}{6-\epsilon}}\|p\|_{L^{\frac{6-\epsilon}{6-2\alpha}}(\mathcal{C}(\frac R2, R))}\|\vu\|_{L^{\frac{6-\epsilon}{3-\alpha}}(\mathcal{C}(\frac R2, R))}.
\end{eqnarray*}
If $1+\frac{\epsilon}{3}< \alpha< \frac{5}{3}+\frac{2}{9}\epsilon$, the power of the parameter $R$ is then negative and we have $\underset{R\to+\infty}{\lim}I_c=0$, while if $\alpha=\frac{5}{3}+\frac{2}{9}\epsilon$ we use the fact that  $\|p\|_{L^{\frac{6-\epsilon}{6-2\alpha}}(\mathcal{C}(\frac R2, R))},\|\vu\|_{L^{\frac{6+\epsilon}{3-\alpha}}(\mathcal{C}(\frac R2, R))}\underset{R\to +\infty}{\longrightarrow 0}$. Now in the case $\alpha=1+\frac{\epsilon}{3}$, we have $\vu\in L^3(\R)$, $p\in L^{\frac{3}{2}}(\R)$and we can proceed as in (\ref{EstimateR3}). We thus have $\underset{R\to+\infty}{\lim}I_c=0$.\\
\item[$\bullet$] If $\frac{3}{5}<\alpha<1$ we have $\vu\in L^{\frac{6+\epsilon}{3-\alpha}}(\R)$ and by Lemma \ref{Lemme_InformationPression} we also have $p\in  L^{\frac{6+\epsilon}{6-2\alpha}}(\R)$. As above, if $1-\frac{\epsilon}{3}< \alpha\leq \frac{5}{3}-\frac{2}{9}\epsilon$, by the H\"older inequality with $\frac{3\alpha-3-\epsilon}{6-\epsilon}+\frac{6-2\alpha}{6-\epsilon}+\frac{3-\alpha}{6-\epsilon}=1$, we obtain
$$I_c\leq C R^{-1+3\frac{3\alpha-3-\epsilon}{6-\epsilon}}\|p\|_{L^{\frac{6-\epsilon}{6-2\alpha}}(\mathcal{C}(\frac R2, R))}\|\vu\|_{L^{\frac{6-\epsilon}{3-\alpha}}(\mathcal{C}(\frac R2, R))}.$$
Note that if $1-\frac{\epsilon}{3}< \alpha< \frac{5}{3}-\frac{2}{9}\epsilon$, the power of the parameter $R$ is negative, while if $\alpha=\frac{5}{3}-\frac{2}{9}\epsilon$ we have $-1+3\frac{3\alpha-3+\epsilon}{6+\epsilon}=0$, but we have $\|\vu\|_{L^{\frac{6+\epsilon}{3-\alpha}}(\mathcal{C}(\frac R2, R))}\underset{R\to +\infty}{\longrightarrow 0}$. Remark also that if $\alpha=1-\frac{\epsilon}{3}$, then $\vu \in L^{\frac{6+\epsilon}{3-\alpha}}(\R)=L^{3}(\R)<+\infty$ by hypothesis, $p\in L^{\frac{3}{2}}(\R)<+\infty$ by the Lemma \ref{Lemme_InformationPression} and we can  proceed as in (\ref{EstimateR3}). In any case we have $\underset{R\to+\infty}{\lim}I_c=0$.
\end{itemize}
\end{itemize}
We have proven that
$$\underset{R\to+\infty}{\lim}I_a=0, \quad \underset{R\to+\infty}{\lim}I_b=0 \quad\mbox{and}\quad \underset{R\to+\infty}{\lim}I_c=0,$$
thus making $R\to +\infty$ in the both sides of the inequality (\ref{Inegalite}) we easily obtain that 
$$\|\vu\|_{\dot{H}^{\frac{\alpha}{2}}}=0,$$
from which we deduce by the Sobolev embeddings that $\|\vu\|_{L^{\frac{6}{3-\alpha}}}=0$ and we finally obtain that $\vu\equiv 0$. Theorem \ref{Theo_Liouville} is proven.  \hfill $\blacksquare$

\mysection{Proof of Theorem \ref{Theoreme_Regularite}}\label{Secc_ProofTheo3}
\begin{itemize}
\item We start proving the first point of the Theorem \ref{Theoreme_Regularite}. Recall that in this case we have $\frac{5}{3}<\alpha<2$. Thus, applying the Leray projector $\mathbb{P}$ to the fractional Navier-Stokes equation and using the divergence free condition, we have the equation $(-\Delta)^{\frac{\alpha}{2}}\vu=-\mathbb{P}(\dv(\vu\otimes \vu))$ which can be rewritten as 
$$\vu=-\frac{\mathbb{P}(\dv(\vu\otimes \vu))}{(-\Delta)^{\frac{\alpha}{2}}}.$$
Now, for some index $\sigma>0$ that will be defined later, we write
$$\|(-\Delta)^{\frac{\sigma}{2}}\vu\|_{L^2}=\left\|(-\Delta)^{\frac{\sigma}{2}}\frac{\mathbb{P}(\dv(\vu\otimes \vu))}{(-\Delta)^{\frac{\alpha}{2}}}\right\|_{L^2}\leq C\|(-\Delta)^{\frac{\sigma-\alpha+1}{2}}(\vu\otimes \vu)\|_{L^2},$$
where we used the boundedness properties of the Leray projector in the $L^2$ space. At this point we apply the product law given in the Lemma \ref{ProductRule} to obtain (since $\vu\in \dot{H}^{\frac{\alpha}{2}}(\R)$)
\begin{equation}\label{Application_LoiProduit}
\|(-\Delta)^{\frac{\sigma-\alpha+1}{2}}(\vu\otimes \vu)\|_{L^2}=\|\vu\otimes \vu\|_{\dot{H}^{\sigma-\alpha+1}}\leq C\|\vu\|_{\dot{H}^{\frac{\alpha}{2}}}\|\vu\|_{\dot{H}^{\frac{\alpha}{2}}}<+\infty,
\end{equation}
as long as $\sigma-\alpha+1=\alpha-\frac{3}{2}$, from which we deduce that $\sigma=2\alpha-\frac{5}{2}$. Now, since $\alpha>\frac{5}{3}$ we have that $\sigma>\frac{\alpha}{2}$. We have thus proved that 
$$\|(-\Delta)^{\frac{\sigma}{2}}\vu\|_{L^2}=\|\vu\|_{\dot{H}^\sigma}<+\infty,$$
which is a gain of regularity. By iterating this process we easily obtain that the solutions of the equation (\ref{NSStationnaire2}) are smooth. 
\item We study now the second point of  Theorem \ref{Theoreme_Regularite}, where we have $1<\alpha\leq \frac{5}{3}$. In this case, we have $\vu\in \dot{H}^{\frac{\alpha}{2}}(\R)$ which seems to be not enough to obtain a gain of regularity when applying the Lemma \ref{ProductRule} in the estimate (\ref{Application_LoiProduit}). To circumvent this issue, we will use the additional hypothesis given by $\vu\in L^\infty(\R)$ and instead of Lemma \ref{ProductRule} we use the Leibniz fractional inequality given in Lemma \ref{FracLeibniz} and in (\ref{Application_LoiProduit}) we write:
$$\|(-\Delta)^{\frac{\sigma-\alpha+1}{2}}(\vu\otimes \vu)\|_{L^2}\leq  C\|(-\Delta)^{\frac{\sigma-\alpha+1}{2}}\vu\|_{L^2}\|\vu\|_{L^\infty},
$$
which is a finite quantity as long as $\sigma-\alpha+1=\frac{\alpha}{2}$, which gives $\sigma=\frac{3}{2}\alpha-1$. Since $1<\alpha\leq \frac{5}{3}$ we have $\sigma>\frac{\alpha}{2}$ and we have obtained a gain of regularity as we have proved that $\vu\in \dot{H}^{\sigma}(\R)$. Again, by iteration we obtain that the solutions of the equation (\ref{NSStationnaire2}) are smooth. 
\end{itemize}
Theorem \ref{Theoreme_Regularite} is now proven. \hfill $\blacksquare$


\end{document}